\documentclass[a4paper,12pt]{amsart}
\usepackage[margin=3cm]{geometry}
\usepackage{url}
\newtheorem{theorem}{Theorem}[section]

\newtheorem{proposition} [theorem]{Proposition}
\newtheorem{lemma}[theorem]{Lemma}
\newtheorem{definition}[theorem]{Definition}
\newtheorem{corollary}[theorem]{Corollary}
\newcommand{\dquer}{\overline\partial}

\newcommand{\nm}{\nabla_{min}}
\newcommand{\nms}{\nabla_{min}^*}
\newcommand{\nmm}{\nabla_{max}}
\newcommand{\nmms}{\nabla_{max}^*}

\numberwithin{equation}{section}
\begin{document}
\title{An idea on proving weighted Sobolev embeddings.}

\author{Klaus Gansberger}

\thanks{Supported by the FWF-grant P19667.}

 \address{K. Gansberger, Institut f\"ur Mathematik, Universit\"at Wien,
Nordbergstrasse 15, A-1090 Wien, Austria}
\email{klaus.gansberger@univie.ac.at}

\maketitle

\begin{abstract} ~\\
This article contains a characterization of when certain weighted Sobolev spaces on $\Bbb R^n$ embed compactly into $L^2(\mathbb R^n , \varphi )$. This characterization is in terms of derivatives of the weight function $\varphi$ and involves the Wiener capacity, as it is obtained from reformulating the problem in terms of resolvent properties of Schr\"odinger operators. This reformulation also works for general domains. 
\end{abstract}

\section{Introduction and Results.}~\\

The aim of the present paper is to characterize when certain weighted Sobolev spaces embed compactly into the weighted Lebesgue space
$$L^2(\mathbb R^n , \varphi )=\{ f:\mathbb R^n \to \mathbb C \ : \ \int_{\mathbb R^n}|f|^2\, e^{-\varphi}\,d\lambda < \infty \},$$
where $\varphi$ is a weight function that is typically smooth and bounded from below. For bounded domains $\Omega$ with smooth boundary $\partial\Omega$, one has the Rellich -- Kondrachov Theorem assuring that the classical unweighted Sobolev space of order one is compactly embedded in $L^2(\Omega)$. This has many consequences --- let us point out some of them, to illustrate the importance of this Theorem in the various different fields of mathematics, but this is no means meant to be complete. Compact Sobolev embedding Theorems play a role  in the theory of non-linear partial differential equations for dealing with certain boundary value problems. Theorems of this type can also be used to show discreteness of the spectrum of strictly elliptic differential operators. The reader can find details on such applications for instance in \cite{ee} and \cite{et}. Moreover these Embedding Theorems are important in statistics for showing the existence of an orthonormal set of Nonlinear Principal Components, see e.g. \cite{sal} and its references. Coming from complex analysis, they classically are used to prove the equivalence of compactness in the $\dquer$-Neumann problem to the existence of so-called compactness estimates, see for instance \cite{str}. \\
Since the Rellich -- Kondrachov Theorem enjoys many applications, there was a great attempt to generalize it in various different directions, for instance replacing the smoothness condition by certain weaker geometric ones, e.g. the cone condition. If one tries to consider unbounded domains, one quickly realizes by taking transverses of a function that one can not hope for a compact embedding unless the domain is getting sufficiently thin at infinity. To overcome this problem, one can either change the Definition of the Lebesgue space $L^p$ by allowing more general convex functions than $x\mapsto |x|^p$ generating the so-called Orlicz-spaces, or one can introduce weight functions. In this paper, we will investigate the second possibility further.\\
So let us cosider a smooth weight function $\varphi$ and define the following notions of weighted Sobolev spaces.
 
\begin{definition}
\label{sob1}
For a domain $\Omega$ and a weight function $\varphi$ define 
$$H^1 (\Omega, \varphi):= \{ f\in L^2(\Omega, \varphi) \ : \ \frac{\partial f}{\partial x_j} \in 
 L^2(\Omega, \varphi) \ {\text{for}} \ 1 \le j\le n \}.$$
 Here, we understand the derivatives in the sense of distributions and equip the space with the norm
 $$\|f\|_{1,\varphi}^2=\| f\|^2_\varphi+ \sum_{j=1}^n\left\|\frac{\partial f}{\partial x_j}\right\|_\varphi ^2.$$
Let moreover $H^1_0(\Omega, \varphi)$ be the closure of $\mathcal C_0^\infty (\Omega)$, i.e., the space of smooth functions with compact support in $\Omega,$ under the norm defined above.
\end{definition}

Note that any two weight functions that induce equivalent norms on $L^2(\Omega , \varphi )$ will define the same space. So it is not much of a restriction to just consider smooth weights, as modifying a given weight in a bounded way will not change the spaces.\\
We will also be interested in the following closely related Definition of a weighted Sobolev space. This notion appeared for the first time in \cite{bdh}, motivated by a question by F. Mignot. Note that the vector fields $X_j$ appearing in the Definition are the formal adjoints of $-\frac{\partial}{\partial x_j}$ in $L^2(\Omega, \varphi)$.
 
\begin{definition}
\label{sob2}
For $j=1,\dots ,n$ let $X_j=\frac{\partial }{\partial x_j} - \frac{\partial \varphi}{\partial x_j}$ and set
$$H^1(\Omega, \varphi,\nabla\varphi):=\{ f\in L^2(\Omega, \varphi) \ :\ X_jf \in  L^2(\Omega, \varphi)\ {\text{for all}} \ 1 \le j\le n\},$$
with the norm
$$\|f\|^2_{ \varphi ,\nabla \varphi}= \|f\|^2_{ \varphi }+\sum_{j=1}^n \|X_jf\|^2_\varphi .$$
Let again $H^1_0(\Omega, \varphi,\nabla\varphi)$ be the closure of $\mathcal C_0^\infty (\Omega)$ under the respective norm.
\end{definition}

Compact injections of this kind of Sobolev spaces have for instance been used in complex analysis to show compactness results for the $\dquer$-Neumann problem on unbounded domains, see \cite{gh} and \cite{ga}.\\
Note that both Definitions \ref{sob1} and \ref{sob2} extend the classical notion of a Sobolev space in the sense that they define the same space when $\Omega$ is bounded and $\varphi\in\mathcal C^\infty(\overline\Omega)$. Note also that the norm in $H^1(\Omega, \varphi,\nabla\varphi)$ is related to but not the same as taking $\|fe^{-\varphi/2}\|_1$.\\
In case that $\varphi$ is a subharmonic weight function, i.e., $\triangle\varphi\ge0$, it holds that $H_0^1(\Omega, \varphi,\nabla\varphi)\subset H^1_0(\Omega, \varphi)$, see \cite{ga}, Lemma 4.3. But in general, there is no such relation for the spaces $H^1(\Omega, \varphi,\nabla\varphi)$ and $H^1(\Omega, \varphi)$.
\vskip 0.5cm

The main idea of this article is the following reformulation in terms of Schr\"odinger operators and its consequences. Surprisingly, although most of the ideas used in this paper are known in spectral theory and have been used at least implicitly, this point of view nevertheless seems to be new and there is no literature about it. We combine classical ideas that in fact go back to Witten with a compactness criterion from \cite{g}. This allows us to express compactness of the injection of Sobolev spaces by resolvent properties of certain Schr\"odinger operators. The resolvent of an operator $A$ is the map 
$$\mathfrak R(z, A)= (A-z)^{-1},$$
which is defined for all $z\in \Bbb C\setminus \sigma(A)$, where $\sigma (A)$ is the spectrum of $A$. An operator is said to have compact resolvent, if $(A-z)^{-1}$ is compact for one (or equivalently for all) $z\in\Bbb C\setminus\sigma(A)$.\\
At least in the case of $\Bbb R^n$, compactness of the resolvent of a Schr\"odinger operator is understood --- a fact that will yield a characterization in terms of the weight function. \\
Our reduction for $\Bbb R^n$ is the following.

\begin{proposition}
\label{equiv}
The injection  $H^1 (\Bbb R^n, \varphi)\hookrightarrow L^2(\Bbb R^n, \varphi)$ is compact if and only if the Schr\"odinger operator $\mathcal S_1=-\triangle+V_1$ has compact resolvent, where the potential $V_1$ equals $V_1=\frac{1}{4}|\nabla\varphi|^2-\frac{1}{2}\triangle\varphi$.\\
The injection  $H^1 (\Bbb R^n, \varphi,\nabla\varphi)\hookrightarrow L^2(\Bbb R^n, \varphi)$ is compact if and only if $\mathcal S_2=-\triangle+V_2$ has compact resolvent, where $V_2=\frac{1}{4}|\nabla\varphi|^2+\frac{1}{2}\triangle\varphi$.\\
\end{proposition}

In fact more is true: $\mathcal S_1$ is not only a Schr\"odinger operator, but a Witten Laplacian, see Section \ref{connection} for details.\\
In the sequel we will give several conclusions from Proposition \ref{equiv}, obtained by applying known criteria for compactness of the resolvent of non-magnetic Schr\"odinger operators. If the potential $V$ of a Schr\"odinger operator is semi-bounded from below, i.e., $V(x)\ge -C$ for some $C>0$, it is a result of A. Mol\v canov \cite{mol} that it is possible to express discreteness of the spectrum in terms of the Wiener capacity. This result has been refined in \cite{kms}. See Section \ref{connection} for the Definitions of capacity and the Mol\v canov functional $M_\gamma$. 

\begin{theorem}
\label{char}
Let $Q_d$ denote a cube with edges of length $d$ that are parallel to the coordinate axes and suppose that $V_1\ge-C$ for some $C>0$. Then there is $c>0$ such that the injection 
$$H^1 (\Omega, \varphi)\hookrightarrow L^2(\Omega, \varphi)$$
is compact if and only if there is $d_0>0$ such that for all $d\in (0,d_0)$
$$d^{-n}M_\gamma(Q_d,V_1)\to\infty\quad\text{as}\quad Q_d\to\infty,$$
where $\gamma=cd^2/g(d)$, $g:(0,d_0)\to(0,\infty)$ is any function with $g(d)\to0$ for $d\to0$ and $d^2\le g(d)$ and finally $Q_d\to \infty$ means that the center of the cube goes to infinity.\\
The analogous statement holds for the injection $H^1 (\Bbb R^n, \varphi,\nabla\varphi)\hookrightarrow L^2(\Bbb R^n, \varphi)$ after replacing $V_1$ by $V_2$.
\end{theorem}

Let us point out that in Mol\v canov's result, there is no need to restrict oneself to cubes. One can equivalently take any other system of sets defining the Euclidean topology on $\Bbb R^n$.\\
Interestingly, as the properties of the embedding in Theorem \ref{char} are invariant under equivalent weights, this is also true for the limit of the Mol\v canov functional.
\vskip 0.5cm

Although it is natural that it arises, conditions involving capacity are in practice hard to handle. It was proven in \cite{ks}, Section 6.1, that it is possible to replace the capacity by the Lebesgue measure $\lambda$ in order to get sufficient conditions:

\begin{corollary}
\label{mes}
Let $r>0$ and $\Bbb B(x,r)$ denote the ball with center $x$ and radius $r$. Suppose that $V_1$ is semi-bounded below and suppose that there is $\gamma>0$ such that for any discrete sequence $(x_l)_l\subset\Bbb R^n$ and any compact sets $F_l\subset \Bbb B(x_l,r)$ with $\lambda(F_l)\le\gamma r^n$ it holds that 
$$\int_{\Bbb B(x_l,r)\setminus F_l}\left(\frac{1}{4}|\nabla\varphi|^2-\frac{1}{2}\triangle\varphi\right)d\lambda(y)\to\infty\quad\text{as}\quad k\to\infty.$$
Then the injection $H^1 (\Bbb R^n, \varphi)\hookrightarrow L^2(\Bbb R^n, \varphi)$ is compact.\\
The analog statement holds for the case of $H^1 (\Bbb R^n, \varphi,\nabla\varphi)$.
\end{corollary}

By semi-boundedness of the potential, we can without loss of generality assume that it is positive. Thus we increase the domain of integration in Theorem \ref{char} to get the following obvious necessary condition.

\begin{corollary}
\label{int}
Suppose that $V_1\ge-C$ for some $C>0$. If the injection $H^1 (\Omega, \varphi)\hookrightarrow L^2(\Omega, \varphi)$ is compact, then for any $r>0$ it holds that
$$\int_{\Bbb B(x,r)}\left(\frac{1}{4}|\nabla\varphi|^2-\frac{1}{2}\triangle\varphi\right)d\lambda\to\infty\quad\text{as}\quad |x|\to\infty.$$
If $V_2\ge-C$ for some $C>0$ and the injection $H^1 (\Bbb R^n, \varphi,\nabla\varphi)\hookrightarrow L^2(\Bbb R^n, \varphi)$ is compact, then for any $r>0$ 
$$\int_{\Bbb B(x,r)}\left(\frac{1}{4}|\nabla\varphi|^2+\frac{1}{2}\triangle\varphi\right)d\lambda\to\infty\quad\text{as}\quad |x|\to\infty.$$
\end{corollary}

{\bf Remark. }The conditions say that the mean value of $V_1$ or $V_2$ in balls going to infinity is necessary for a compact embedding. If $V_2(x)\to \infty$ for $|x|\to\infty$, it was shown in \cite{bdh} that this is sufficient for a compact embedding of $H^1 (\Bbb R^n, \varphi,\nabla\varphi)$. The analog statement for $H^1 (\Bbb R^n, \varphi)$ was essentially shown in \cite{jo}, Proposition 6.2, by using the same idea (apart from a minor flaw in the proof that can be fixed).\\
\vskip 1cm

\section{The connection to a Schr\"odinger operator.}~\\
\label{connection}

In this section, we reduce the question of compact injections to resolvent properties of certain Schr\"odinger operators.  We combine ideas which were used in \cite{g} with classical ones that go back to E. Witten, by showing that the spaces of our interest can be interpreted as the domain of certain singular operators on $L^2(\Bbb R^n, \varphi)$. These operators are unitarily equivalent to certain Schr\"odinger operators, so-called Witten Laplacians. For a good introduction to the theory of Witten Laplacians we refer to \cite{he}.\\
Let us first recall the following Proposition, which reformulates a part of Theorem 3 from \cite{g} to our setting. 

\begin{proposition}
\label{compact}
Let $T$ be a linear partial differental operator with smooth coefficients acting on $dom(T)\subset L^2(\Omega, \varphi)$. Suppose furthermore that $T$ is closed and densely defined. Let $T_\varphi^*$ be the adjoint of $T$ in $L^2(\Omega, \varphi)$ and set $P=T^*_\varphi T$.\\
Then the follwing are equivalent:
\begin{enumerate}
\item $P$ has compact resolvent.
\item The injection $j_\varphi$ of the space $dom(T)$ equipped with the inner product $\langle u,v\rangle_T=\langle Tu,Tv\rangle_\varphi$ into $ L^2(\Omega, \varphi)$ is compact.
\end{enumerate}
\end{proposition}

For completeness, let us give the short argument which is based on an idea of E. Straube that appeared in \cite{str}. Both compactness of the resolvent and of the injection imply $\dim\ker(T)<\infty$ and in particular closedness of the range of $T$, so it suffices to show that $dom(T)\cap(\ker(T))^\perp\hookrightarrow L^2(\Omega, \varphi)$ is compact. $T^*_\varphi T |_{(\ker(T))^\perp}$ is continuously invertible, so we can without loss of generality assume that $P=T^*_\varphi T$ is invertible and call its inverse $P^{-1}$. Now for all $u,v\in dom(T)$ it holds
$$\langle u,v\rangle_\varphi=\langle u,j_\varphi v\rangle_\varphi=\langle j_\varphi^* u,v\rangle_T,$$
while on the other hand
$$\langle u,v\rangle_\varphi=\langle PP^{-1}u,v\rangle_\varphi=\langle TP^{-1}u,Tv\rangle_\varphi=\langle P^{-1}u,v\rangle_T.$$
Hence, $P^{-1}=j_\varphi^*$ as an operator to $dom(T)$ and consequently  $P^{-1}=j_\varphi\circ j_\varphi^*$ as an operator to $L^2(\Omega, \varphi)$, which shows the equivalence of $(1)$ and $(2)$ in the Proposition.\\

This allows us to reformulate the problem in Proposition \ref{equiv} in terms of singular differential operators on $L^2(\Bbb R^n, \varphi)$. 

\begin{lemma}
\label{lemma1}
The following equivalences hold:
\begin{enumerate}
\item The injection $H^1 (\Bbb R^n, \varphi)\hookrightarrow L^2(\Bbb R^n, \varphi)$ is compact if and only if 
$$P_1= -\triangle+\sum_{j=1}^n\frac{\partial\varphi}{\partial x_j}\frac{\partial}{\partial x_j}$$ 
acting on $L^2(\Bbb R^n, \varphi)$ has compact resolvent.
\item The injection $H^1 (\Bbb R^n, \varphi,\nabla\varphi)\hookrightarrow L^2(\Bbb R^n, \varphi)$ is compact if and only if 
$$P_2= -\triangle+\sum_{j=1}^n\frac{\partial\varphi}{\partial x_j}\frac{\partial}{\partial x_j}+\triangle\varphi$$ 
acting on $L^2(\Bbb R^n, \varphi)$ has compact resolvent.
\end{enumerate}
\end{lemma}

{\em Proof. } Let $\nabla : L^2(\Bbb R^n, \varphi)\to \oplus_{j=1}^n L^2(\Bbb R^n, \varphi)$ be given by $\nabla f=(\frac{\partial f}{\partial x_1},\dots ,\frac{\partial f}{\partial x_n})$, where we think of $\nabla$ as the maximal extension of the operator initially defined on $\mathcal C_0^\infty(\Bbb R^n)$. Then by Definition $dom (\nabla)=H^1 (\Bbb R^n, \varphi)$ and $\dim\ker(\nabla)\le1$. In particular the range of $\nabla$ is closed. So by Proposition \ref{compact}, the injection $H^1 (\Bbb R^n, \varphi)\hookrightarrow L^2(\Bbb R^n, \varphi)$ is compact if and only if $\nabla^*_\varphi \nabla $ has compact resolvent.\\
A standard computation shows that $\nabla_\varphi^* : \oplus_{j=1}^n L^2(\Bbb R^n, \varphi)\to L^2(\Bbb R^n, \varphi)$ is given by $\nabla _\varphi^* (u_1,\dots,u_n)=-\sum_{j=1}^n\left(\frac{\partial u_j}{\partial x_j}-\frac{\partial \varphi}{\partial x_j}\right)$. Let us emphasize that we use at this point density of $\mathcal C_0^\infty(\Bbb R^n)$ in the graph norm to be able to do integration by parts. Thus 
$$\nabla^*_\varphi \nabla f= -\triangle f+\sum_{j=1}^n\frac{\partial\varphi}{\partial x_j}\frac{\partial f}{\partial x_j},$$
which shows $(1)$.\\
One can prove $(2)$ verbatim after defining  $T: L^2(\Bbb R^n, \varphi)\to \oplus_{j=1}^n L^2(\Bbb R^n, \varphi)$,  $T f=(\frac{\partial f}{\partial x_1}-\frac{\partial \varphi}{\partial x_1}f,\dots ,\frac{\partial f}{\partial x_n}-\frac{\partial \varphi}{\partial x_n}f)$
\begin{flushright}
$\square$
\end{flushright}

{\bf Remark. } Looking at the injection $H^1 (\Bbb R^n, \varphi)\hookrightarrow L^2(\Bbb R^n, \varphi)$, it is obvious that compactness is invariant under equivalent weights. The same is true for the injection $H^1 (\Bbb R^n, \varphi,\nabla\varphi)\hookrightarrow L^2(\Bbb R^n, \varphi)$, as the following argument shows: Let $G_\varphi: \oplus_{j=1}^n L^2(\Bbb R^n, \varphi)\to L^2(\Bbb R^n, \varphi)$ be the canonical solution operator (i.e., the one mapping to the orthogonal complement of the kernel) to the equation $T^*f=g$ in $L^2(\Bbb R^n, \varphi)$, where $T$ is from above. Then it is easy to see that the resolvent of $T^*T$ is compact if and only if $G_\varphi$ is compact. Now if $\psi$ is an equivalent weight, the inclusion $\iota: L^2(\Bbb R^n, \psi)\to L^2(\Bbb R^n, \varphi)$ is continous. Since the equation $T^*f=g$ does not depend on the weight, $G_\psi=\iota\circ G_\varphi\circ\iota^{-1}$ defines a compact solution operator to $T^*$ on $L^2(\Bbb R^n, \psi)$, which implies compactness of the injection $H^1 (\Bbb R^n, \psi,\nabla \psi)\hookrightarrow L^2(\Bbb R^n, \psi)$.
\vskip 0.5cm

The computation in the following Lemma dates back to Witten. But let us include it for completness.

\begin{lemma}
\label{lemma2}
$P_1$ acting on $L^2(\Bbb R^n, \varphi)$ is unitarily equivalent to a non-magnetic Schr\"odinger operator $\mathcal S_1= -\triangle+V_1$ acting on $L^2(\Bbb R^n)$, more precisely it holds 
$$e^{-\varphi/2}P_1e^{\varphi/2}= -\triangle+V_1,$$
where $V_1=\frac{1}{4}|\nabla\varphi|^2-\frac{1}{2}\triangle\varphi$. Moreover we have
$$e^{-\varphi/2}P_2e^{\varphi/2}= -\triangle+V_2,$$
with $\mathcal S_2= -\triangle+V_2$ and $V_2=\frac{1}{4}|\nabla\varphi|^2+\frac{1}{2}\triangle\varphi$.\\
\end{lemma}

{\em Proof. } The proof is a straight forward computation. Let us first compute
\begin{align*}
e^{-\varphi/2}\triangle e^{\varphi/2}=&e^{-\varphi/2}\nabla\left(\frac{1}{2}\nabla\varphi e^{\varphi/2}+e^{\varphi/2}\nabla\right)\\
=&e^{-\varphi/2}\left(\frac{1}{2}\triangle\varphi e^{\varphi/2}+\frac{1}{4}|\nabla\varphi|^2 e^{\varphi/2}+\frac{1}{2}\nabla\varphi e^{\varphi/2}\nabla+\frac{1}{2}\nabla\varphi e^{\varphi/2}+e^{\varphi/2}\triangle\right)\\
=&\triangle +\nabla\varphi \nabla+\frac{1}{2}\triangle\varphi +\frac{1}{4}|\nabla\varphi|^2 .
\end{align*}
Moreover we have
\begin{align*}
e^{-\varphi/2}\left(\nabla\varphi\nabla\right) e^{\varphi/2}=&\nabla\varphi\nabla+\frac{1}{2}|\nabla\varphi|^2,
\end{align*}
which proves the Lemma.
\begin{flushright}
$\square$
\end{flushright}

 {\em Proof of Proposition \ref{equiv}.} The spectra of unitarily equivalent operators coincide, so $P_1$ has compact resolvent, i.e., empty essential spectrum, if and only if $\mathcal S_1$ has. This follows by Lemma \ref{lemma2} and similarly for $P_2$. Thus Lemma \ref{lemma1} completes the proof of the Proposition.
\begin{flushright}
$\square$
\end{flushright}

{\em Proof of Theorem \ref{char}.} To proof the Theorem, it suffices to cite Mol\v canov's result. We do this in the generalized version of \cite{kms}. To this end we first need to define some more notions.\\
Let $\Omega\subset\mathbb R^n$ be open, and let $F\subset \Omega$ be a compact set. The capacity of $F$ with respect to $\Omega$ is
$$\operatorname{cap}_\Omega(F)=\inf \left\lbrace\int_{\mathbb R^n} |\nabla u(x)|^2 d\lambda : u\in Lip_0(\Omega),\  u\equiv 1\ \text{on}\ F\right\rbrace,$$
where $Lip_0(\Omega)$ is the space of all Lipschitz functions with compact support in $\Omega$ and $\lambda$ denotes the Lebesgue measure. By $Q_d$, we always denote a cube with sidelength $d$ and edges parallel to the coordinate axes. If we drop the subscript, $\operatorname{cap}(F)$ is the capacity of $F$ with respect to $\Bbb R^n$ if $n\ge3$, or with respect to $Q^\circ_{2d}$ if $n=2$, where $Q_d$ is the smallest square containing the compact set $F$ and the center of $Q_{2d}$ is the center of $Q_d$. Note that in the Definition of the capacity one could equivalently take the infimum over $u\in\mathcal  C_0^\infty(\Omega),\ 0\le u\le1$, see \cite{ks}, Section 3.\\
The Mol\v canov functional is defined by
\begin{equation}
M_\gamma(Q_d,V)=\inf _{F\subset Q_d}\left \lbrace \int_{Q_d\setminus F}V(x)d\lambda(x) : \operatorname{cap}(F)\le\gamma\operatorname{cap}(Q_d)\right\rbrace,
\end{equation}
where $0<\gamma<1$. Due to properties of the capacity the infimum will not change if we restrict it to compact sets $F$ which are the closures of smooth open subsets of $Q_d$. Now applying Mol\v canov's result in the more general version of Kondratiev, Maz$^\prime$ya and Shubin (see \cite{kms},  Theorem 1.2 or Theorem \ref{kms} below) we get Theorem \ref{char} directly from Proposition \ref{equiv}.

\begin{theorem}
\label{kms}
Let $\mathcal S=-\Delta+V(x)$ and suppose that $V\in L^1_{loc}$ is semi-bounded from below. Then there is $c>0$ such that $\mathcal S$ has compact resolvent if and only if there is $d_0>0$ such that for all $d\in (0,d_0)$
$$d^{-n}M_\gamma(Q_d,V)\to\infty\quad\text{as}\quad Q_d\to\infty,$$
where $\gamma=cd^2/g(d)$ and $g:(0,d_0)\to(0,\infty)$ is any function with $g(d)\to0$ for $d\to0$ and $d^2\le g(d).$ 
\end{theorem}
\vskip 0.5cm

{\em Proof of Corollary \ref{mes}.} The Corollary follows from Proposition \ref{equiv} and Theorem 6.1 in \cite{ks}.
\begin{flushright}
$\square$
\end{flushright}

{\em Proof of Corollary \ref{int}.} If the $V_j$'s are semibounded from below, we can without loss of generality assume that they are positive since this will not change discreteness of the spectrum. Thus Corollary \ref{int} follows from Corollary \ref{char} by increasing the domain of integration.
\begin{flushright}
$\square$
\end{flushright}
\vskip 0.5cm

{\bf Example.} Suppose that the weight $\varphi$ is a polynomial. Then the injection $H^1 (\Omega, \varphi)\hookrightarrow L^2(\Omega, \varphi)$ is compact if and only if 
\begin{equation}
\label{nabla}
\int_{\Bbb B(x,1)}|\nabla\varphi|^2d\lambda\to \infty
\end{equation}
for $|x|\to \infty$.\\
First we observe that condition \eqref{nabla} is equivalent to
\begin{equation}
\label{nabla2}
\int_{\Bbb B(x,1)}\left(\frac{1}{4}|\nabla\varphi|^2-\frac{1}{2}\triangle\varphi\right)d\lambda\to \infty,
\end{equation}
since the first term is the one with the higher degree. By the same reason is $\mathcal S=-\triangle+\frac{1}{4}|\nabla\varphi|^2-\frac{1}{2}\triangle\varphi$ always bounded from below. Although we can not directly use Corollary \ref{int}, it follows from general spectral theory that compactness of the resolvent $\mathcal S$ implies $\langle \mathcal S\chi_j,\chi_j\rangle\to\infty$ for each normed sequence going weakly to zero and belonging to the domain of $\mathcal S$. So taking cut-off functions over balls that have uniformly bounded second order derivatives makes \eqref{nabla2} as a necessary condition for a compact injection immediate.\\
On the other hand, \eqref{nabla2} is sufficient for compactness of the resolvent, as follows for instance from the Fefferman-Phong Lemma as it was given in \cite{aube}.
\vskip 1cm

\section{General domains.}~\\

For completeness, let us also consider general unbounded domains, where the situation is more subtle. In the case of $\Bbb R^n$, there is no difference between the spaces $H^1 (\Bbb R^n, \varphi)$ and $H^1_0 (\Bbb R^n, \varphi)$ and similarly for $H^1 (\Bbb R^n, \varphi,\nabla\varphi)$. On general domains, we must distinguish between those spaces and be more careful when applying our idea, since the $X_j$'s from Definition \ref{sob2} are only the formal adjoints of $\partial /\partial x_j$. Moreover, there are different choices of a closed extension of $\nabla$ and we gain not as much from our reformulation, as few criteria for compactness of the resolvent of a Schr\"odinger operator on domains with boundary are known.\\
We restrict ourselves to the case of $H^1 (\Omega, \varphi)$ and note that $H^1 (\Omega, \varphi, \nabla\varphi)$ can be treaten analogously.\\

\begin{lemma}
Let $\nm$ be the minimal closed extension of $\nabla$. Then $P=\nm^*\nm$ coincides with the Friedrich's extension of the operator defined by the quadratic form
$$Q^D_\varphi (f,f)=\int_\Omega|\nabla f|^2\ e^{-\varphi}\ d\lambda, \quad f\in \mathcal C_0^\infty(\Omega).$$
\end{lemma}

{\em Proof.} On $\mathcal C_0^\infty(\Omega)$, the operator defined by $Q^D_\varphi$ clearly coincides with $\nms\nm$. If $\nm$ is the minimal extension of $\nabla$, then $\nms$ is maximal and we have $dom (\nm)=H_0^1(\Omega,\varphi)$ and $dom(\nms)=\{f\in \oplus_{j=1}^nL^2(\Omega,\varphi)\ | \ \nms f\in L^2(\Omega,\varphi) \}$. The domain of $P$ is $dom(P)=\{ f\in dom (\nm)\ |\ \nm f\in dom(\nms)\}$, thus $dom (P)=H_0^1(\Omega,\varphi)\cap\{ f\ |\ \nm f\in dom(\nms)\}$.\\
Now by general facts about  Friedrich's extension, the domain of the extension is the intersection of the form domain  with the domain of the adjoint of the operator one wants to extend, see e.g. \cite{weid}. Consequently, this is $H_0^1(\Omega,\varphi)\cap\{ f\in H_0^1(\Omega,\varphi)\ |\ \nm f\in dom(\nms)\}$, which proves the Lemma.
\begin{flushright}
$\square$
\end{flushright}

Analogously  one shows:
\begin{lemma}
Let $\nmm$ be the maximal closed extension of $\nabla$. Then $P=\nmms\nmm$ coincides with the Friedrich's extension of the operator defined by the quadratic form
$$Q^N_\varphi (f,f)=\int_\Omega|\nabla f|^2\ e^{-\varphi}\ d\lambda, \quad f\in \mathcal C^\infty(\overline\Omega).$$
\end{lemma}

This allows us to argue similarly to the case of $\Bbb R^n$. 

\begin{proposition}
The injection $H^1_0 (\Omega, \varphi)\hookrightarrow L^2(\Omega, \varphi)$ is compact if and only if the Dirichlet realization of the Schr\"odinger operator $\mathcal S_1=-\triangle+V_1$ has compact resolvent, where $V_1$ is from Proposition \ref{char}.\\
The injection $H^1 (\Omega, \varphi)\hookrightarrow L^2(\Omega, \varphi)$ is compact if and only if the von Neumann realization of $\mathcal S_1$ has compact resolvent.
\end{proposition}

{\em Proof.} Let us show the first statement. The domain of $\nm$ is $H^1_0 (\Omega, \varphi)$, so the injection is compact if and only if $P=\nms\nm$ has compact resolvent. The Friedrich's extension of $Q^D_\varphi$ corresponds to the Dirichlet realization of $P$, see for instance \cite{he2}, Chapter 2. The latter one in fact is unitarily equivalent to the Dirichtlet realization of $\mathcal S_1$, as follows from the discussions in \cite{he}, Chapter 2.5. 
\begin{flushright}
$\square$
\end{flushright}

{\bf Remark. } The reader can find a characterization of when the Dirichlet realization of a Schr\"odinger operator in a domain $\Omega$ has compact resolvent in \cite{kms2}.

\end{document}